\documentclass[a4paper,11pt]{article}
\usepackage{amsmath,amscd,amsfonts,amssymb,epsf,latexsym}

\newtheorem{theorem}{Theorem}

\newtheorem{lemma}[theorem]{Lemma}
\newtheorem{corollary}[theorem]{Corollary}
\newtheorem{remark}[theorem]{Remark}
\newtheorem{definition}[theorem]{Definition}
\newtheorem{example}[theorem]{Example}

\def\proof{\noindent{\bf Proof.\ }}

\def\qed{~\hbox{$\Box$}\s}
\def\D{{\mathcal D}}
\def\C{{\mathbb C}}
\def\Q{{\mathbb Q}}
\def\F{{\mathbb F}}
\def\K{{\mathbb K}}

\def\Hc{{\cal H}}

\def\s{\vskip6pt}

\begin{document}
\baselineskip20pt

\title{\bf Purity at the end}

\author{ Andrzej Weber\thanks{Research supported by a KBN grant
NN201 387034.}
\\
\small Department of Mathematics of Warsaw University\\
\small Banacha 2, 02-097 Warszawa, Poland\\
\small aweber@mimuw.edu.pl}

%\subjclass{***}
\date{September 2009}

\maketitle

\begin{abstract}  \baselineskip 18pt
We consider smooth completion of algebraic
manifolds. Ha\-ving some information about its singular
completions or about completions of its images we prove purity of
cohohomology of the set at infinity. We deduce also some
topological properties. The work is based on the study of perverse
direct images for algebraic maps.

\vskip6pt\noindent{\itshape Key words and phrases.} Algebraic
varieties, weight filtration, perverse sheaves.

\end{abstract}

\section{Introduction}
This paper is inspired by the important paper of De Cataldo and
Miglio\-ri\-ni \cite{dCM}, where the authors give a proof of
Decomposition Theorem of \cite[Th\'eor\`eme 6.2.5]{BBD} using
Hodge-theoretic methods. As a side result they obtain some
contractibility criteria for subvarieties. We present a
generalizations of these criteria. Our proof is based on Purity
Theorem of Gabber \cite{G} or stability of pure perverse sheaves
with respect to intermediate extension \cite[Corollarie
5.3.2]{BBD}: for a pure perverse sheaf $K$ supported by an open
set the intermediate extension $j_{!*}K$ is pure. We use formally
the properties of the derived category of sheaves and our argument
works equally well in the setup of Weil sheaves (as in \cite{BBD})
and for mixed Hodge modules of \cite{Sa}.

We study open smooth algebraic varieties and their smooth
completions. Let $X=\overline U$ be a completion of $U$.
 Our reasoning is localized around a connected
component $Z$ of the set $X\setminus U$, which we call an and of
$U$. Knowing that some  image $V=f(U)$ of $U$ admits a (possibly
singular) completion $Y$ we study the weight filtration in the
cohomology of the link of $Z$.

Let us explain the case when $\overline f:X\to Y$ is a resolution
of an isolated singularity, $Z\subset X$ is the exceptional locus,
$U=X\setminus  Z$, $V=Y\setminus \overline f(Z)$, $f=\overline
f_{|U}$. We have an exact sequence
$$\begin{matrix}
&&_\alpha&&_\beta&&_\gamma&\cr
 \to&H^{k-1}(U)&\to& H^k(X,U)&\to&H^k(X)&\to&H^{k-1}(U)&\to\,.\cr\end{matrix}$$
The following argument is valid for complex varieties, but with
some suitable changes\footnote{Although the notion of the link
$L_W(Y)$ of a subvariety $W$ in $Y$ is ambiguous (and has no
meaning for varieties over $\F_q$) it is justified to talk about
the cohomology $H^*(L_W(Y);\K):=H^*(W;i^*j_*j^*\K_Y)$, where
$j:Y\setminus W\to Y$ and $i:W\to Y$ are inclusions.} can be
performed in the category of varieties over a finite
characteristic field (see the proof of Theorem \ref{main1}). We
present it to give a topological motivation of further
constructions. We replace $Y$ by a conical neighbourhood of the
singular point $\{p\}=f(Z)$. Then $X$, the resolution of $Y$, is a
manifold with boundary, which retracts to the exceptional set $Z$.
The boundary $\partial X$ is homeomorphic to the link of $p$ in
$Y$, denoted by $L_p(Y)$. The open set $U$ retracts to $\partial
X\simeq L_p(Y)$. Let $m=\dim X$. By Poincar\'e duality we have
$H^k(X,U)\simeq H_{2m-k}(Z)$.
 The long exact sequence
can be rewritten in the following way:
$$\begin{matrix}
&&_\alpha&&_\beta&&_\gamma&\cr
 \to &H^{k-1}(L_p(Y))&\to &
H_{2m-k}(Z)&\to &H^k(Z)&\to &H^{k}(L_p(Y))&\to .
 \end{matrix}$$
 All the cohomology groups appearing in the diagram are
equipped with  weight filtration. Since $Z$ is complete
\begin{itemize}
 \item $H^k(X,U)$ is of weight $\geq k$, (the isomorphism
 with $H_{2m-k}(Z)=\left(H^{2m-k}(Z)\right)^*$ shifts the weight by $2m$),
 \item $H^k(Z)$ is of weight $\leq k$,
\end{itemize}
Moreover, by purity of the intersection sheaf $IC_Y$ we have
\begin{itemize}
 \item for $k<m$ the group $H^k(L_p(Y))\simeq\Hc^ k(IC_Y[-m])_p$ is of weight
 $\leq k$,
 \item by duality, for $k\geq m$ the group $H^k(L_p(Y))\simeq
 H^{k+1}_{\{p\}}(Y;IC_Y[-m])$ is of weight $ \geq k$,
\end{itemize}
see \cite{DM}.  It follows that for $k\leq \dim Y$ the map
$\alpha$ vanishes and the map $\beta$ is injective. (Dually, for
$k\geq \dim Y$ the map $\gamma$ vanishes and $\beta$ is
surjective.) The case of isolated singularities was already
studied in \cite{GM} as a corollary from the decomposition
theorem. Additionally, from injectivity of $\alpha$ we obtain that
$H^k(X,U)$ is pure of weight $k$ for $k\leq \dim X$. By duality we
have

\begin{theorem} \label{ga} For a resolution of an isolated singularity
$\overline f:X\to Y$ the cohomology of the exceptional set
$H^k(Z)$ is pure of weight $k$ for $k\geq \dim X$.\end{theorem}

We would like to emphasize that by  \cite{W} the condition that
$H^k(Z)$ is pure of weight $k$ depends only on topology of $Z$.
\bigskip

The Theorem \ref{ga} is  proved by de Cataldo and Migliorini,
 as a corollary from their proof of
Decomposition Theorem. In fact it is enough to assume that the map
$f$ is semismall, \cite[Th. 2.1.11]{dCM}. We generalize the
Theorem \ref{ga} without any assumption for the map $f:U\to V$.
Instead we have a condition for the degree $k$ for which  the
cohomology group $H^k(Z)$ is pure. The Decomposition Theorem of
\cite{BBD} or \cite{Sa} is involved in our proof but we do not
rely on it, although we have to  use equally strong arguments as
purity of the intermediate extension of a pure perverse sheaf.

Our argument is formal. We list below the formal properties of
weights which we use in the proof. These properties hold equally
well for Weil sheave of \cite{BBD} and for mixed Hodge modules of
M.~Saito \cite{Sa}. In the first case we apply $\ell$-adic \'etale
cohomology  for varieties defined over a finite field $\F_q$
($\ell\not=char(\F_q)$), that is
$H^*(X_0;F_0):=H^*_{et}(X_0\otimes_{\F_q}\overline{\F}_q;F)$,
where $F_0\in\D(X_0)$ is a Weil sheaf, see the notation of
\cite[\S5.1]{BBD}. These cohomology groups are vector spaces over
$\overline\Q_\ell$, the closure of the field of $\ell$-adic
numbers. In the second case the cohomology is a vector space over
rationals or complex numbers. Depending on the context let $\K$
denote $\overline\Q_\ell$, $\Q$ or $\C$.
\bigskip

Notation:

\begin{description}

\item{--} $\D(X)$ denotes the category of Weil sheaves or the
category of mixed Hodge modules. The objects of $\D(X)$ are called
sheaves.

\item{--} We drop the letter $R$ for the direct image of an object
of $\D(X)$, i.e.~we write $f_*F$ instead of $Rf_*F$ for $f:X\to
Y$.

\item{--} On the other hand for the constant sheaf $\K_X$ we write
$Rf_*\K_X$ since $\K_X$ is a sheaf in the usual sense.

\item{--} For the truncation functor we write $\tau_{<d}$ instead
of $\tau_{\leq d-1}$. \end{description}

We list the formal properties of weights which are used in the
proofs.

\begin{description}

\item{(0.1)} functors $f_!$ and $f^*$ preserve the weight
condition $\leq w$,

\item{(0.2)} in particular if $F\in \D(X)$ is of weight $\leq w$,
then $H^k_c(X;F)$ is of weight $\leq w+k$,

\item{(0.3)} functors $f_*$ and $f^!$ preserve the weight
condition $\geq w$,

\item{(0.4)} in particular if $F\in \D(X)$ is of weight $\geq w$,
then $H^k(X;F)$ is of weight $\geq w+k$,

\item{(0.5)} Verdier duality switches the weight conditions $\geq
w$ and $\leq -w$,

\item{(0.6)}  Verdier duality exchanges $!$ with $*$, i.e. $D\circ
f^!=f^*\circ D$ and $D\circ f_!=f_*\circ D$,

\end{description}

\section{Main Theorem}

The set of components $U_\infty=\pi_0( X\setminus  U)$ for a
normal completion $X$ of a normal algebraic variety $U$ does not
depend on the completion. An element of this set is called an end
of $U$.
 A  map of algebraic varieties  which is proper induces a map of their
ends.
 A completion of an end $\eta\in U_\infty$ is the component
 $Z\subset X\setminus  U$ corresponding to $\eta$. We say that the
completion $Z$ of $\eta$ is smoothing if $X$ is smooth in a
neighbourhood of $Z$. Below we give the exact statement of our
main theorem.

\begin{theorem} \label{teo2} Suppose we have a proper map of
varieties $f:U\to V$ with smooth $U$ of dimension $m$  and let
$r=r(f)$ be the defect of semismallness. Let $\eta$ be an end of
$U$. Suppose that $f(\eta)\in V_\infty$ admits a completion with
the variety at the end of dimension $d$. Then for any smoothing
completion $Z$ of $\eta$ the cohomology $H^k(Z)$ is pure for the
degrees $k\geq m+r+d$. Moreover, for that range of degrees the
restriction map $H^k(X,U)\to H^k(Z)$ is surjective.
\end{theorem}

The notion of the defect of semismallness was invented in
\cite{dCM}, we recall it in Definition \ref{dss}. If $f$ is a
fibration then $r(f)$ is the dimension of the fiber.

\bigskip In the proof we reduce the general situation to the case
when $f:U\to V$ extends to $\overline f:X\to Y$. Then Theorem
\ref{teo2} follows from Corollary \ref{main1c} which is the dual
version of Theorem \ref{main1}. There we deal with an arbitrary
pure sheaf $F$ on $Y$ instead of the sheaf $R\overline f_*\K_X$
and instead
 of extracting geometric conditions of $f$ we record vanishing range of
perverse cohomology $^p{\cal H}^k(F)$.

\section{Defect of semismallness}
Assume that the variety $U$ is smooth. The following discussion is
to clarify the  relation between the decomposition
\cite[Th\'eor\`eme 6.2.5]{BBD}
$$Rf_*\K_U=\bigoplus IC_{S_\alpha}(L_\alpha)[d_\alpha]$$
and the geometric properties of the map $f$. We recall the
definition of a number which measures to what extend $Rf_*\K$ is
not perverse.

 Set
$$V^i=\{y\in V:\dim f^{-1}(y)=i\}\,.$$

\begin{definition} \rm \cite[Def 4.7.1]{dCM} \label{dss} The {\it defect of
semismallness} of a map $f:U\to V$ is the integer
$$r(f)=\max_{i: V^i\not=\emptyset}\{2i+\dim V^i-\dim U\}\,.$$
\end{definition}
Note that
 \begin{itemize}
 \item  $r(f)=0$ if and only if $f$ is semismall.
 \item If $f$ is a
fibration, then $r(f)$ is equal to the dimension of the fiber.
\end{itemize}
The defect of semismallness controls the number of perverse
derived images which appear in the Decomposition Theorem. To see
that we need to recall the  perverse t-structure in the derived
category $\D(V)$, according to which for a closed smooth variety
$W\subset Y$ the sheaf $\K_W[\dim W]$ is perverse. Precisely:
 $$F\in{^p\D(V)}^{\leq 0}\quad{\rm if}\;\forall\, s\; \dim({\rm
Supp}\,{\cal H}^{s}F)\leq -s\,.$$ Then
$$F\in{^p\D(V)}^{\leq r}\quad{\rm if}\;\forall\, s\; \dim({\rm
Supp}\,{\cal H}^{s}F)\leq r-s\,.$$
 If we apply this condition to $F=Rf_*\K_U[m]$ we obtain
 $$Rf_*\K_U[m]\in{^p\D(V)}^{\leq r}\quad{\rm if}\;\forall\,t\; \dim({\rm
Supp}\,R^t f_*\K_U)\leq m+r-t\,.\eqno{(*)}$$
 Note that for $V^i$ introduced before the definition of the
 defect of semismallness we have
$$\overline{V^i}={\rm Supp}\,R^{2i}f_*\K_U\,.$$
Hence the condition (*) is satisfied for even $t=2i$ provided that
$\dim V^i\leq m+r-2i$. For odd $t=2i-1$ we have ${\rm
Supp}\,R^{2i-1}f_*\K_U\subset {\rm Supp}\,R^{2i}f_*\K_U$, so (*)
follows from the condition for $t=2i$. Therefore the condition (*)
is satisfied for all $t$ provided that the map has the defect of
semismallnes at most to $r$.

 We conclude that for any proper map
$f:U\to V$ with $r(f)=r$ we have a decomposition
$$Rf_*\K_U[m]\simeq \bigoplus_{s=-r}^r {^p{\cal H}}^s( Rf_*\K_U[m])[-s]\,.$$
Note that the perverse cohomology for degrees smaller  than $-r$
vanish since $Rf_*\K_U[m]$ is self-dual
(i.e.~$D(Rf_*\K_U[m])\simeq Rf_*\K_U[m]$).

\section{Proof of the main theorem}

Let $Y$ be an algebraic variety. Let $W\subset Y$ be a closed
subvariety of dimension $d$ and let $V=Y\setminus W$. Denote by
$j$ the inclusion $V\hookrightarrow Y$.

\begin{lemma} \label{per}
Let $P$ be a perverse sheaf on $Y\setminus  W$, which is of weight
$\leq w$, then for $y\in W$ the cohomology stalks
$\Hc^k((j_*P)_y)$ are of weight $\leq k+w$ for $k<-d$. In another
words $\tau_{< -d}\,j_*P$ is of weight $\leq w$.
\end{lemma}

\proof The natural map
$$\tau_{< -d}\,j_{!*}P\to \tau_{<-d}\,j_*P$$
is an isomorphism by construction \cite[Proposition 2.1.17]{BBD}.
The sheaf $j_{!*}P$ is of weight $\leq w$ by the stability of
weight with respect to the intermediate extension \cite[Corollarie
5.3.2]{BBD}. \qed

We fix a natural number $r$ and consider the sheaves which have
the perverse cohomology bounded below by $-r$. We obtain the
corollary:

\begin{corollary}\label{per2} Let $F\in {^p\D^{\geq -r}(V)}$, be a
sheaf which is of weight $\leq w$. Then for $y\in W$ the
cohomology stalks $\Hc^k((j_*F)_y)$ are of weight $\leq k+w$ for
$k<-d-r$. In another words $\tau_{< -d-r}\,j_*P$ is of weight
$\leq w$.
\end{corollary}

\proof If $F$ is of weight $\leq w$ then its preverse cohomology
$^p\Hc^k(F)$ is of weight $\leq w+k$. In the framework of Weil
sheaves this is \cite[Th\'eor\`eme 5.4.1]{BBD}. In Saito
construction this is the definition \cite[4.5.1]{Sa}. By usual
spectral sequence argument we deduce the claim from Lemma
\ref{per}.\qed

\begin{remark}\rm We will apply Corollary \ref{per2} for pure sheaves $F$. In
that case we do not have to use spectral sequences since $F$ is
isomorphic to the sum of its perverse cohomology, i.e.
$$F=\bigoplus_{s\geq r} {^p\Hc^s}(F)[-s]\,,$$
by \cite[Th\'eor\`eme 5.4.5]{BBD} or by \cite[formula 4.5.4]{Sa}.
\end{remark}

Let $i:W\hookrightarrow Y$ be the inclusion. Suppose that $W$ is
complete.

\begin{lemma} \label{per3} Let $F\in {^p\D^{\geq -r}(V)}$
be a  sheaf
 on $V$, which is of
weight $\leq w$. Then the cohomology $H^k(W;i^*j_*F)$ is of weight
$\leq w+k$ for $k<-d-r$.\end{lemma}

\proof The map $H^k(W;i^*\tau_{<-d-r}\,j_*F)\to H^k(W;i^*j_*F)$ is
an isomorphism for $k<-d-r$. The conclusion follows from (0.1),
(0.2) and Corollary \ref{per2} .\qed

Here is our key statement:

\begin{theorem} \label{main1}
Let $G\in \D(Y)$ such that $j^*G\in{^p\D^{\geq -r}(Y)}$. Assume
that $W$ is complete. If $G$ is pure of weight $w$, then
$H^k(W;i^!G)$ is pure
 of weight $w+k$ for $k\leq -d-r$. Moreover the natural map
 $H^k(W;i^!G)\to H^k(W;i^*G)$ is injective.\end{theorem}

 \proof
 We have a distinguished triangle
$$\begin{matrix} i_*i^!G&\longrightarrow&\phantom{gg} G\cr
      _{[+1]}\hfill    \nwarrow&&\swarrow\hfill&\cr
          & j_*j^*G.\cr\end{matrix}$$
Let us restrict this triangle to $W$. We have $i^*i_*i^!G=i^!G$
and we obtain the the triangle
$$\begin{matrix} i^!G&\longrightarrow&i^*G\cr
      _{[+1]}\hfill    \nwarrow&&\swarrow\hfill&\cr
          & i^*j_*j^*G.\cr\end{matrix}$$
Consider the associated sequence of cohomology {\small
$$\begin{matrix}&_\alpha& &_\beta& &_\gamma\cr
 H^{k-1}(W;i^*j_*j^*G)&\to& H^k(W;i^!G)&\to& H^k(W;i^*G)&\to&
H^{k}(W;i^*j_*j^*G)\,\cr\end{matrix}$$} (compare the exact
sequence from the Introduction). By (0.3) the sheaf $i^!G$ is of
weight $\geq w$. By (0.4) the
 cohomology group $H^k(W;i^!G)$ is of weight $\geq w+k$.
On the other hand $i^*G$ is of weight $\leq w$ by (0.1) and  since
$W$ is complete the
 cohomology group $H^k(W;i^*G)=H^k_c(W;i^*G)$ is of weight $\leq w+k$.
By Lemma \ref{per3} applied to $F=j^*G$  the cohomology
$H^{k-1}(W;i^*j_*j^*G)$ is of weight $\leq w+k-1$ for $k-1\leq
-d-r$. It follows that the map $\alpha$ is trivial. Therefore the
map $\beta$ is injective and the group $H^k(W;i^!G)$ is pure.\qed

\begin{corollary} \label{main1c} Suppose $G\in \D(Y)$ is pure of weight $w$
and $j^*G\in{^p\D^{\leq r}(Y)}$. Then the cohomology
$H^k(W;i^*G)$ is pure for $k>d+r$ and the map
$\gamma:H^k(W;i^!G)\to H^k(W;i^*G)$ is surjective.\end{corollary}

\proof We apply the Theorem \ref{main1} for the dual sheaf $DG$.
By the property (0.5) of Verdier duality the assumptions Theorem
\ref{main1} are satisfied. Hence by (0.6) and Theorem \ref{main1}
$$H^k(W;i^!DG)=H^{-k}(W;i^*DG)^*,\quad
H^k(W;i^*DG)=H^{-k}(W;i^!DG)^* $$  is pure for $k\geq d+r$.\qed

We specialize Theorem \ref{main1} to the following situation.
 Let $\overline f:X\to Y$ be a proper map  and let $W\subset Y$ be a
complete subvariety of dimension $d$.  Denote by $Z$ the inverse
image $f^{-1}(W)$:
 $$\begin{matrix}
 & _{j'} &&_{ i'}\cr
 U&\hookrightarrow & X & \hookleftarrow &Z \cr
 &\phantom{.}\cr
\phantom{^f} \downarrow^f&&\phantom{^f}\downarrow^{\overline
 f}&&\phantom{^{f^\infty}}\downarrow^{f^\infty}\cr
  &\phantom{.}\cr
 V&\hookrightarrow&Y &\hookleftarrow  &W&.\cr
 &^j& &^i\cr\end{matrix}$$
We assume that $U$ is smooth. Let $r=r(f)$ be the defect of
semismallness of $f=\overline f_{|U}$.
 The sheaf $G=\overline f_*IC_X$ is pure of weight $m=\dim
X$ and it satisfies the assumption of the Theorem \ref{main1} and
Corollary \ref{main1c}. We have $f^\infty_*i'^!=i^!\overline f_*$
and $f^\infty_*i'^*=i^*\overline f_*$.
 We obtain the following:

\begin{theorem} \label{main2} With the  notation as above we have
 \begin{enumerate}
 \item $H^k(Z;i^!IC_X)$ is pure of weight $m+k$ for $k\leq -d-r$,
 \item $H^k(Z;i^*IC_X)$ is pure of weight $m+k$ for $k\geq d+r$,
 \item the map $\beta:H^k(Z;i^!IC_X)\to H^k(Z;i^*IC_X)$
 is injective for $k\leq -d-r$ and surjective for $k\geq d+r$.
 \end{enumerate}
 \end{theorem}

When $X$ is smooth of dimension $m$ we translate the conclusion of
Theorem \ref{main2} to the cohomology with coefficients in the
constant sheaf $\K_Z$. We obtain:
 \begin{enumerate}
 \item $H^k_Z(X;\K)=H^k(X,U;\K)$
 is pure of weight $k$ for $k\leq m-d-r$,
 \item $H^k(Z;\K)$ is pure of weight $k$ for $k\geq m+d+r$,
 \item the map $\beta:H^k(X,U;\K)\to H^k(Z;\K)$
 is injective for $k\leq m-d-r$ and surjective for $k\geq m+d+r$.
 \end{enumerate}

We have obtained the conclusion of Theorem \ref{teo2} except that
we have to get rid of the assumption that the map $f:U\to V$
extends to $\overline f:X\to Y$. To this end consider the variety
$$\widehat X=\overline{graph(f)}\subset X\times Y\,.$$
We apply Theorem \ref{main2} for $\widehat f=pr_Y:\widehat X\to Y$
and deduce purity of $H^*(\widehat Z;i^!IC_{\widehat X})$
 Let
$g=pr_X:\widehat X\to X$ be the projection.
 The composition of the natural maps{\small
 $$\K_X[m]\to g_*\K_{\widehat X}[m]\to g_*IC_{\widehat X}\to g_*D\K_{\widehat
 X}[m]=Dg_*\K_{\widehat X}[m]\to D\K_X[m]=\K_X[m]$$}
 \vskip-10pt
 \noindent is the identity.
 Hence $\K_X[m]$ is a direct summand of $g_*IC_{\widetilde X}$ in a natural way.
 Therefore $H^k(X,U;\K)$ is pure of weight $k$
for $k\leq m-d-r$. The purity of $H^k( Z;\K)$ follows by duality.
Also triviality of the map from
$\alpha:H^{k-1}(Z;i^*j_*j^*\K_X[m])\to H^k(Z;i^!\K_X[m]) $ follows
by comparison with the covering maps for $\widehat X$. We deduce
injectivity  of the map $\beta:H^k(X,U;\K)\to H^k(Z;\K)$ for
$k\leq m-d-r$. The surjectivity for the complementary degrees
follows by duality. This completes the proof of the Theorem
\ref{teo2}.

\section{Fibration case}

It is interesting to analyze the case when $f$ is a fibration. If
the fiber is of dimension $r$, then the defect of semismallness is
equal to $r$. Set $n=\dim V=m-r$. Then $V$ is a smooth variety of
dimension $n$. The derived direct image decomposes
$$Rf_*\K_U=\bigoplus_{s=0}^{2r}R^sf_*\K_U[-s]\,.$$
We have a decomposition of the cohomology of link of $Z$:
$$H^k(L_Z(X))=\bigoplus_{s=0}^{2r}H^{k-s}(L_W(Y);R^sf_*\K_U)\,.$$
By \cite{DM} the cohomolgy of the link $H^k(L_W(Y))$ is of weight
$\leq k$ for $k<m-d$. Lemma \ref{per3} generalized this statement
to the case of twisted coefficients (note that $R^sf_*\K_U$ is
pure of weight $s$). Again we see that $H^k(L_Z(X))$ is of weight
$\leq k$ for $k<m-d=n-r-d$.

\section{Corollaries}

We derive some easy corollaries. Fist of all the restriction map
$H^k(X,U)\to H^k(Z)$ factors through $H^k(X)$, therefore we have:

\begin{corollary} With the notation of Theorem \ref{teo2} for
$k \geq m+r+d$ the restriction map $H^k(X)\to H^k(Z)$ is
surjective.\end{corollary}

In particular when $f:U\simeq V$ we obtain an easy criterion for
contracti\-bi\-li\-ty:

\begin{corollary} Let $Z$ be a complete subset of an algebraic
manifold $X$. Suppose that $rk\,H^k(X)<rk\,H^k(Z)$ then there does
not exist an algebraic variety $Y$ with a subvariety $W\subset Y$
such that $Y\setminus W\simeq X\setminus Z$ and $\dim W\leq
k-m$.\end{corollary}

When $X$ is complete we obtain an exact sequences of pure
cohomology groups:

\begin{corollary} Suppose that $X$ is complete (and smooth as before).
With the notation of Theorem \ref{teo2} for $k > m+r+d$ we have an
exact sequences of pure cohomology groups of weight $k$
 $$0\to H^k(X,Z)\to H^k(X) \to H^k(Z)\to 0\,.$$
 Dually, for $k<m-r-d$ we have an exact sequence of pure
 cohomology groups of weight $k$
 $$0\to H^k(X,U)\to H^k(X) \to H^k(U)\to 0\,.$$
\end{corollary}

Information about the weight structure of $H^*(Z)$ can exclude
some of contractions regardless of the ambient space $X$:

\begin{example}\rm  Let $C={\bf P}^1/\{0,\infty\}$ be the nodal
curve. Then $H^1(C)=\K$ is of weight 0. Let $M$ be a smooth
variety of dimension $k$. Then $H^{2k+1}(M\times C)$ is not pure.
Therefore $Z=M\times C$ cannot be shrunk to a point in a variety
$X$ of the dimension $n\leq 2k+1$. (E.g.~$k=1$, $n=3$.) Another
example of this type is given in \cite[Remark
2.1.13]{dCM}.\end{example}

Note that the condition that $H^k(Z)$ is pure depends only on the
topology of $Z$. According to \cite{W} the bottom  part of the
weight filtration $W_{k-1}H^k(Z)$ is equal to the kernel of the
natural map to intersection cohomology$H^k(Z)\to IH^k(Z)$, and
both invariants depend only on topology. The assumptions of
Theorem \ref{teo2} imply that $H^k(Z)\to IH^k(Z)$ is injective for
$k\geq m+r+d$.

\begin{remark} \rm The case when $f$ is an isomorphism,
$d=0$, i.e.~the case of an isolated singularity was studied by
many authors (e.g. \cite{GM}). Let $Y$ be a conical neighbourhood
of an isolated singularity and $X$ be its resolution, $Z$ the
exceptional set (as in the introduction). The middle part of the
exact sequence of the pair $(X,
\partial X)$ has the form:
$$\begin{matrix}&&_0&&_\simeq&&_0\cr
\to& H^{n-1}(\partial X)&\to&H^n(X,\partial X)
 &\to&H^n(X)&\to&H^n(\partial X)&\to\,.\cr\end{matrix}$$
The isomorphism in the center can be translated to the
nondegenerate intersection form on $H^n(X,U)$. For the degrees
$k<n$ the map $H^k(X)\to H^k(\partial X)$ is a surjection. The
kernel is isomorphic to $H^k(X,U)$. Note that since this kernel is
pure we obtain that $H^k(\partial X)$ is pure if and only if
$H^k(Z)$ is pure. We generalize this remark.
\end{remark}

\begin{corollary} With the assumption of Theorem \ref{teo2}
for $k<m-r-d$ the cohomology of the link $H^k(L_Z(X))$ is pure of
weight $k$ if and only if $H^k(Z)$ is pure of weight
$k$.\end{corollary}

\begin{remark}\rm When $X$ is complete and when
there exists a map $\overline f:X\to Y$ extending $f$ then the
statement of Theorem \ref{teo2} can be obtained directly  applying
Decomposition Theorem. We have
$$R\overline f_*\K_X[m]\simeq \bigoplus_{s=-r}^{r}{^p{\Hc}}^s(\overline f_*\K_{X}[m])
[-s]\oplus F\,,$$
 where $F$ is a sheaf supported by $W$.
 In another words
 $$R\overline f_*\K_X\simeq A\oplus B\,,$$
where $A\in{^p\D}^{\leq r+m}(Y)$ and $B=F[-m]$ is supported by
$W$.
 Counting the dimensions of the supports of ${\Hc}^k(A)$ we find that
$H^k(W;A)=0$ for $k\geq m+r+d$ and
$$H^k(Z)\simeq H^k(W;(\overline f_*\K_{X})_{|W})\simeq H^k(W;B_{|W})\simeq H^k(Y;B)\,.$$
 But the last group is a direct summand in
$H^k(X)=H^k(Y;R\overline f_*\K_X)$. Hence it is pure. By the
exactly the same argument one proves Corollary \ref{main1c}, but
we have chosen the way which allows to trace directly how the
weights of link cohomology is involved.
\end{remark}

\end{document}